\documentclass{elsart}

\usepackage{latexsym}
\usepackage[psamsfonts]{amssymb}
\usepackage{amsfonts}
\usepackage{graphicx}
\usepackage{amsmath,cite}

\journal{eeeeee}

\usepackage{amssymb}

\begin{document}

\newtheorem{tm}{Theorem}[section]
\newtheorem{pp}{Proposition}[section]
\newtheorem{lm}{Lemma}[section]
\newtheorem{df}{Definition}[section]
\newtheorem{tl}{Corollary}[section]
\newtheorem{re}{Remark}[section]
\newtheorem{eap}{Example}[section]

\newcommand{\pof}{\noindent {\bf Proof} }
\newcommand{\ep}{$\quad \Box$}

\newcommand{\al}{\alpha}
\newcommand{\be}{\beta}
\newcommand{\var}{\varepsilon}
\newcommand{\la}{\lambda}
\newcommand{\de}{\delta}
\newcommand{\st}{\stackrel}

\allowdisplaybreaks

\begin{frontmatter}



\title{Some notes on characterizations of compact sets in fuzzy number spaces
}

\author{Huan Huang $^{a}$\corauthref{cor}}
\author{Congxin Wu $^{b}$}
\author{}\ead{hhuangjy@126.com (H. Huang), wucongxin@hit.edu.cn (C. Wu)
}
\address{$^{a}$Department of Mathematics, Jimei
University, Xiamen 361021, China \newline $^{b}$Department of
Mathematics, Harbin Institute of Technology, Harbin 150001, China  }
\corauth[cor]{Corresponding author.}

\date{}

\begin{abstract}
  In this paper, we presents a characterization of compact subsets of the fuzzy number space equipped with the level convergence topology. Based on this, it is shown that compactness is equivalent to sequential compactness on the fuzzy number space equipped with the level convergence topology. Diamond and Kloeden gave a characterization of compact sets in  fuzzy number spaces equipped with the supremum metric,
  Fang and Xue also gave a characterization of compact sets in one-dimensional fuzzy number spaces equipped with supremum metric. The latter characterization is just the one-dimensional case of the former characterization. There exists conflict between the characterization given by us and the characterizations given by the above mentioned authors.
  We point out the characterizations gave by them is incorrect by a counterexample.
\end{abstract}


\end{frontmatter}


\section{Introduction}

The convergences on fuzzy number spaces have been extensively discussed by various authors
\cite{da, gom, kae, roe, h1, fn, gm, wuo, h3}.
One of the most important problems is the characterizations of compact subsets.

In this paper, we present a characterization of compact subsets of the fuzzy number space equipped with level convergence topology. Based on this, we show that compactness is equivalent to sequential compactness on the fuzzy number space equipped with level convergence topology.

Diamond and Kloeden \cite{da} gave a characterization of compact sets in fuzzy number spaces equipped with the supremum metric,
Fang and Xue \cite{fe} gave a characterization of compact sets in one-dimensional fuzzy number spaces equipped with the supremum metric. The compactness criteria given by Fang and Xue is a special case of $m=1$ of the compactness criteria given by Diamond and Kloeden. It is obviously that there exists a contradiction between the characterizations of compact sets given by us and the characterizations given in \cite{da, fe}. We point out the characterizations in \cite{da, fe} are incorrect by a counterexample.

\section{Fuzzy number space}

Let $\mathbb{N}$ be the set of all natural numbers, $\mathbb{R}^m$
be $m$ dimension Euclid space, and $F(\mathbb{R}^m)$ represent all
fuzzy subsets on $\mathbb{R}^m$, i.e. functions from $\mathbb{R}^m$
to $[0,1]$. For details, we refer the readers to references
\cite{wu, da}.

For $u\in F(\mathbb{R}^m)$, let $[u]_{\al}$ denote the $\al$-cut of
$u$, i.e.
\[[u]_{\al}=\begin{cases} \{x\in \mathbb{R}^m : u(x)\geq \al \},\
\al\in(0,1],\\ {\rm supp}\, u=\overline{\{x \in
\mathbb{R}^m: u(x)>0\}}, \al=0.
\end{cases}\]
We call $u\in F(\mathbb{R}^m)$ a fuzzy number if $u$ has the
following properties:

(1) $u$ is normal: there exists at least one $x_{0}\in \mathbb{R}^m$
with $u(x_{0})=1$;\\
(2) $u$ is convex: $u(\la x+(1-\la)y)\geq {\rm min} \{u(x),u(y)\}$
for $x,y \in \mathbb{R}^m$ and $\la \in [0,1];$\\
(3) $u$ is upper semi-continuous;\\
(4) $[u]_0$ is a bounded set in $\mathbb{R}^m$.

The set of all fuzzy numbers is denoted by $E^m$.

Suppose that
$K(\mathbb{R}^m)$ is the set of all nonempty compact sets of $\mathbb{R}^m$
and that $K_c(\mathbb{R}^m)$ is the set of all nonempty compact and convex
set of $\mathbb{R}^m$. The following representation theorem is used
widely in the theory of fuzzy numbers.

\begin{pp}\cite{pr} \label{pr}\
Given $u\in E^m,$ then

(1) \ $[u]_1\not= \emptyset$ and $[u]_\la\in K_c(\mathbb{R}^m)$ for all $\la\in [0,1]$;\\
(2) \ $[u]_\la=\bigcap_{\gamma<\lambda}[u]_\gamma$ for all $\la\in (0,1]$;\\
(3) \ $[u]_0=\overline{\bigcup_{\gamma>0}[u]_\gamma}$.

Moreover, if the family of sets $\{v_\al:\al\in [0,1]\}$ satisfy
conditions $(1)$ through $(3)$ then there exists a unique $u\in E^m$
such that $[u]_{\la}=v_\lambda$ for each $\la\in [0,1].$
\end{pp}

Many metrics and topologies on $E^m$ are based on the well-known
Hausdorff metric. The {\rm Hausdorff} metric $H$ on
   $K(\mathbb{R}^m)$ is defined by:
$$H(U,V)=\max\{H^{*}(U,V),\ H^{*}(V,U)\}$$
for arbitrary $U,V\in K(\mathbb{R}^m)$, where
  $$H^{*}(U,V)=\sup\limits_{u\in U}\,d\, (u,V) =\sup\limits_{u\in U}\inf\limits_{v\in
V}d\, (u,v).$$ Obviously, if $[x_1,\, x_2]$ and $[y_1,\, y_2]$ are
bounded closed intervals of $\mathbb{R}$, then
$$H([x_1,\, x_2],\, [y_1,\, y_2])=\max\{ |x_1-y_1|,\ |x_2-y_2|\}.$$

Throughout this paper, we suppose that the metric on $\mathbb{R}^m$
is the Euclidean metric, and the metric on $K(\mathbb{R}^m)$ is the
Hausdorff metric $H$. The Hausdorff metric has the following
properties.

\begin{pp}
\cite{ke,rm}\label{ec}
$(X,d)$ is a metric space, $K(X)$ is the set of all compact set of
$X$. Then

(1)\ $(X,d)$ complete $\Leftrightarrow \ (K(X), H)$ complete;\\
(2)\ $(X,d)$ separable $\Leftrightarrow \ (K(X), H)$
separable;\\
(3)\ $(X,d)$ compact $\Leftrightarrow \ (K(X), H)$ compact.
\end{pp}

In this paper, we consider two types of convergences on fuzzy number spaces.
\begin{itemize}
\item \ Let $u, u_n\in E, n=1,2,\ldots$. If $\lim_{n\to\infty}d_\infty(u_n,u)=0$, then we say $\{u_n\}$ supremum converges to $u$, denoted by $u_n\st {d_\infty} {\rightarrow} u$, where the supremum metric $d_\infty$ is defined by
\[
d_\infty(u,v)=\sup_{\al\in[0,1]}H([u]_\al,\, [v]_\al)
\]
 for
all $u,v\in E^m$.
\item \
Let $u\in E^m$ and let $\{u_\xi: \xi\in D \}$ be a net in
$E^m$, where $D$ is a direct set. If $\lim\limits_{\xi\in
D}H([u_\xi]_\al,\ [u]_\al)=0$
 for each $\al\in[0,1]$, then we say $\{u_\xi\}$ level converges to $u$,
denoted by $\lim\limits_{\xi \in D}u_\xi=u(l)$ or $u_\xi \st l \rightarrow
u. $

\end{itemize}

The supremum metric convergence is stronger than the level convergence on
$E^m$, i.e. if $\{u_n\}$ supremum metric converges to $u$,
then it also level converges to $u$.

We use $(E^m,d_\infty)$ or $(E^m,\tau(l))$ to denote the fuzzy number space $E^m$ equipped with the supremum metric $d_\infty$ or equipped with the topology $\tau(l)$ induced by
level convergence, respectively.

\section{Characterizations of compact sets and sequentially compact sets in $(E^m,\tau(l))$}

We give characterizations of compact sets and sequentially compact sets, respectively, in $(E^m, \tau(l))$. Based on this, we show that compactness is equivalent to sequential compactness on $(E^m, \tau(l))$. We need some propositions and lemmas at first.

\begin{pp}
\cite{h1}\label{h1}
$(E^m, \tau(l))$ is a Hausdorff space and satisfies the first
countability axiom.
\end{pp}

\begin{lm}\label{sc}
Each compact set of $(E^m, \tau(l))$ is sequentially compact.
\end{lm}

\pof \ By Proposition \ref{h1}, $(E^m, \tau(l))$ satisfies the first
countability axiom, from the basic topology, every countable compact
set of $(E^m, \tau(l))$ is sequentially compact. Since a compact set
is obviously countable compact, and thus each compact set of $(E^m,
\tau(l))$ is sequentially compact. \ep

We say that a set $S$ is relatively compact if it has compact closure.

A set $U\subset E^m$ is said to be uniformly support-bounded if there is a compact set $K\subset \mathbb{R}^m$
such that $[u]_0\subset K$
for all $u\in U$.

Let $\mathcal{F}$ be a family of
functions from $S \subset \mathbb{R}$ to $(K_c(\mathbb{R}^m), H)$. Then
\begin{itemize}

\item $\mathcal{F}$ is said to be equi-left-continuous at $\al$ if for each $
\varepsilon>0$ there exists $ \delta(\al, \varepsilon)
>0$ such that $ H(f(\alpha ) , f(\al') ) <\varepsilon $ whenever
$f\in \mathcal{F}$ and $ \al' \in [\al-\delta, \al]$.

\item  $\mathcal{F}$ is said to be equi-right-continuous at $\al$ if for each $
\varepsilon>0$ there exists $ \delta(\al, \varepsilon)
>0$ such that $ H(f(\alpha ) , f(\al') ) <\varepsilon $ whenever
$f\in \mathcal{F}$ and $ \al' \in [\al, \al+\delta]$.
\end{itemize}
We say that $\mathcal{F}$ is equi-left (right)- continuous on $S$ if
it is equi-left (right)- continuous at each point of $S$. Note that $[u]_\bullet$ (where the $\bullet$ may stand for any
subscript) can be seen as functions
from $[0,1]$ to $K_c(\mathbb{R}^m)$.

\begin{lm}\label{ccc}
A subset $U$ of $(E^m, \tau(l))$ is relatively compact if and only if
the following conditions are satisfied:

(1) $U$ is uniformly support-bounded.\\
(2) $\{[u]_\bullet: u\in U\}$ is equi-left-continuous on $(0, 1]$
and equi-right-continuous at 0.
\end{lm}

\pof \ {\sl Necessity.} \ If $U$ is relatively compact in $(E^m,
\tau(l))$, then, by Lemma \ref{sc}, $\overline{U}$ is sequentially
compact in $(E^m, \tau(l))$, and thus $\{[u]_0: u\in \overline{U}
\}$ is compact in $K_c(\mathbb{R}^m)$. So $\{[u]_0: u\in
\overline{U} \}$ is bounded in $K_c(\mathbb{R}^m)$, then obviously
$U$ is uniformly support-bounded, i.e. condition (1) holds.

Now we prove condition (2). In the opposing case where $\{
[u]_{\al}: u\in U\}$ is not equi-left-continuous at $\al_0\in
(0,1]$. Then there exists $\varepsilon_0>0$ and two sequences
$\{u_n\} \subseteq U$ and $\{\al_n\} \subseteq (0,1]$ with
$\al_n\rightarrow \al_0-,\ n=1,2,\ldots$ such that
\begin{equation}\label{fnfv}H([u_n]_{\al_n},   [u_n]_{\al_0})>\varepsilon_0.\end{equation}
Since $\overline{U}$ is compact, by Lemma \ref{sc}, $\overline{U}$
is sequentially compact. We may assume without loss of generality
that $u_n \st{l}{\rightarrow}  u_0\in \overline{U}$. Note that for a
given $\beta<\al_0$, there is an $N$ such that $\al_0>\al_n>\beta$
for all $n>N$, hence $[u_n]_{\al_0}\subseteq [u_n]_{\al_n}\subseteq
[u_n]_\beta$ for all $n>N$, and thus by \eqref{fnfv}
$$H([u_0]_\beta, [u_0]_{\alpha_0})=\lim_m H([u_n]_{\beta}, [u_n]_{\alpha_0})\geq \lim_m H([u_n]_{\al_n},   [u_n]_{\al_0})\geq \varepsilon_0$$
for all $\beta < \al_0 $, this contradicts with
$[u_0]_{\al_0}=\bigcap_{\beta<\al_0}[u_0]_\beta $. Hence
$\{[u]_\bullet: u\in U\}$ is equi-left-continuous on $(0,1]$.
Similarly, we can prove that $\{[u]_\bullet: u\in U\}$ is
equi-right-continuous at 0.

{\sl Sufficiency.} \ Notice that $(E^m, \tau(l))$ can be seen as a
subset of the product space $\prod_{\al\in [0,1]}
(K_c(\mathbb{R}^m), H)$. Let $\overline{U}$ be the closure of $U$
in $\prod_{\al\in [0,1]} (K_c(\mathbb{R}^m), H)$. Given $ v \in
\overline{U }$, there is a net $\{u_\xi: \xi\in D\}$ of $U$ such that
$v= \lim_{\xi\in D}u_\xi$. Then obviously
\begin{equation}[v]_\al \in K_c(\mathbb{R}^m), \ [v]_\mu\subseteq
[v]_\nu \label{c1} \end{equation} for all $\al\in [0,1]$ and $\mu
\geq \nu$. Given $\gamma\in (0,1]$ and $\varepsilon>0$, from the
equi-left-continuity of $\{[u]_\bullet: u\in U \}$ at $\gamma$,
there is a $\delta>0$ such that
$$H([u_\xi]_\gamma, [u_\xi]_{\gamma-\delta})<\varepsilon/3$$ for all
$\xi\in D$. Since $v=\lim_{\xi\in D}u_\xi$, there exists $k\in D$ such
that
$$H([v]_\gamma, [u_k]_\gamma)<\varepsilon/3, \
H([v]_{\gamma-\delta},[u_k]_{\gamma-\delta})<\varepsilon/3.$$ Thus
$$H([v]_{\gamma},[v]_{\gamma-\delta})
\leq H([v]_{\gamma}, [u_k]_{\gamma})+  H([u_k]_\gamma,
[u_k]_{\gamma-\delta} )+ H([u_k]_{\gamma-\delta},
[v]_{\gamma-\delta})< \varepsilon,$$ and so
\begin{equation}\lim_{\delta\to 0}
H([v]_{\gamma},[v]_{\gamma-\delta})=0 \label{c2}
\end{equation} for all $\gamma\in (0,1]$. Combined with \eqref{c1}
and \eqref{c2}, we know
\begin{equation}[v]_\al=\bigcap_{\beta<\al}[v]_\beta \label{c3}\end{equation} for all $\al \in
(0,1]$. Similarly, we can prove that
\begin{equation}[v]_0=\bigcup_{\beta>0}[v]_\beta. \label{c4}\end{equation} Then $v\in E^m$ from Proposition
\ref{pr} and \eqref{c1},\eqref{c3} and \eqref{c4}. So $\overline{U }
\subset E^m$ from the arbitrariness of $v\in \overline{U}$. This
means that the closure of $U$ in $\prod_{\al\in
[0,1]}(K_c(\mathbb{R}^m), H) $ is just the closure of $U$ in
$(E^m, \tau(l))$.

Since $U$ is uniformly support-bounded, then $\{[u]_\al: u \in U \}$
is bounded in $(K_c(\mathbb{R}^m), H)$ for each $\al\in [0,1]$. By
Proposition \ref{ec}, $\overline{\{[u]_\al: u \in U \}}$ is compact
in $(K_c(\mathbb{R}^m), H)$ for each $\al\in [0,1]$, then from the
Tychonoff product theorem $\prod_{\al\in [0,1]} \overline{\{[u]_\al: u \in U
\}} $ is compact in $\prod_{\al\in [0,1]}(K_c(\mathbb{R}^m), H) $.
So $\overline{U}\subset \prod_{\al\in [0,1]}\overline{ \{[u]_\al:
u\in U\} }$ is compact in $\prod_{\al\in [0,1]} (K_c(\mathbb{R}^m),
H)$. Since $\overline{U}\subset E^m$, $\overline{U}$ is also a
compact set in $(E^m,\tau(l))$. \ep

Now, we arrive at one of the main results of this section.

\begin{tm}\label{cc3}
A subset $U$ of $(E^m, \tau(l))$ is compact if and only if the
following conditions are satisfied:

(1) $U$ is closed in $(E^m, \tau(l))$.\\
(2) $U$ is uniformly support-bounded.\\
(3) $\{[u]_\bullet: u\in U\}$ is equi-left-continuous on $(0, 1]$
and equi-right-continuous at 0.
\end{tm}

\pof \ Note that $(E^m, \tau(l))$ is a Hausdorff space, so $U$ is
compact if and only $U$ is closed and relatively compact. The
remainder part of proof follows from Lemma \ref{ccc} immediately. \ep

Fang and Huang \cite{h1} proposed a characterization of compact set
in $(E^m, \tau(l))$. They used concepts ``eventually
equi-left-continuous'' and ``eventually equi-right-continuous''.
\begin{itemize}
\item
 A net $\{u_k\}_{k\in D}$ in $(E^m, \tau(l))$ is said to be eventually equi-left-continuous at $\al\in (0, 1]$,
 if for each $\varepsilon>0$, there exist a $k_0\in D$ and a $\delta>0$ such that $H([u_k ]_{\al-\delta}, [u_k ]_\al)< \varepsilon$ for all
 $k\geq k_0$.
\item
 A net $\{u_k\}_{k\in D }$ in $(E^m, \tau(l))$ is eventually equi-right-continuity at $\al\in [0,1)$, if for each $\varepsilon>0$,
there exists a $k_0\in D$ and a $\delta>0$ such that
$H([u_k]_{\al+\delta}, [u_k]_\al)<\varepsilon$ for all $k\geq k_0$.
\end{itemize}

They \cite{h1} gave the following compact characterization on $(E^m, \tau(l))$.

\begin{pp}\label{chh}
A closed subset $U$ of $(E^m, \tau(l))$ is compact if and only if
 the following conditions are satisfied.

 (1) $U$ is uniformly support-bounded.  \\
 (2) Each net in $U$ has a subnet which is eventually equi-left-continuous on $(0, 1]$ and
 eventually equi-right-continuous at 0.
\end{pp}

The readers may compare
 the condition (3) in Theorem \ref{cc3} with the condition (2) in Proposition \ref{chh}.

\begin{tl}
Suppose that $f$ is a continuous function from $[a,b]$ to $(E^m,
\tau(l))$, then $\{[f(x)]_\bullet: x\in [a,b]\}$ is
equi-left-continuous on $(0,1]$ and equi-right-continuous at $0$.
\end{tl}

\pof \ Since $[a,b]$ is a compact subset of $\mathbb{R}$, we have
$f[a,b]$ is a compact set in $(E^m, \tau(l))$. The desired results
 follows immediately from Theorem \ref{cc3}. \ep

\begin{tm}\label{cs2}
A subset $U$ of $(E^m, \tau(l))$ is sequentially compact if and only
if the following statements are true.

(1) $U$ is closed in $(E^m, \tau(l))$.\\
(2) $U$ is uniformly support-bounded.\\
(3) $\{[u]_\bullet: u\in U\}$ is equi-left-continuous on $(0, 1]$
and equi-right-continuous at 0.
\end{tm}

\pof \ {\sl Necessity.} Given a limit point of $u$ of $U$, since
$(E^m, \tau(l))$ is first countable, there is a sequence $\{u_n,
n=1,2\ldots\}$ of $U$ such that $u=\lim_{n\to\infty}u_n$, and then
 $u\in U$ according to the sequential compactness of $U$. Thus $U$
is a closed set from the arbitrariness of $u$. So statement (1)
holds. Statements (2) and (3) can be proved similarly as in Lemma
\ref{ccc}.

{\sl Sufficiency.} By Theorem \ref{cc3}, if statements (1), (2) and
(3) hold, then $U$ is compact, and thus $U$ is sequentially compact
from Lemma \ref{sc}. \ep

The following statement is another main results of this section.

\begin{tm}
A subset $U$ of $(E^m, \tau(l))$ is compact if and only if it is
sequentially compact.
\end{tm}

\pof \ The desired result follows from Theorems \ref{cc3} and \ref{cs2}. \ep

\section{Characterizations of compact sets in $(E^m, d_\infty)$}

Many authors discussed the characterizations of compact sets in $(E^m, d_\infty)$.
There are many interesting conclusions. However, we find that some of those results is incorrect.

The support function $u^*: [0,1]\times S^{n-1} \to \mathbb{R}$ of $u\in E^m$ is defined by
$$u^*(\alpha,p)=\sup\{   <p,x>: \ x\in [u]_\al   \}.$$

Diamond and Kloeden (Proposition 8.2.1 in \cite{da}) have presented the following compactness criteria of sets in $(E^m, d_\infty)$.

\begin{tm}\label{dao}
A closed set $U$ of $(E^m, d_\infty)$
is compact if and only if

(1)\ $U$ is uniformly support-bounded, and\\
(2)\
$U^*=\{u^*: u\in U\}$ is equi-left-continuous on [0,1] uniformly in $p\in S^{n-1}$, i.e. given $\al\in [0,1]$, for each $\varepsilon>0$, there is a $\delta>0$
such that $ u^*(\alpha,p) \leq u^*(\beta,p) \leq  u^*(\alpha,p)+\varepsilon $ for all $\beta\in [\al-\delta,\alpha]$, $p\in S^{n-1}$, and $u\in U$.
\end{tm}

\begin{re}{\rm
From the properties of support function, we know that condition (2) in Theorem \ref{dao} is equivalent to $\{[u]_\bullet:\ u\in U\}$ are equi-left-continuous on $(0,1]$.
}

\end{re}

Fang and
Xue (Theorem 2.3 of \cite{fe}) gave the following characterization
of compact subsets in $(E^1, d_\infty)$:

\begin{tm}\label{fae}
A subset $U$ in $(E^1, d_\infty)$ is compact if and only if the
following three conditions are satisfied:

(1) $U$ is uniformly support-bounded; \\
(2) $U$ is a closed subset in $(E^1, d_\infty)$;\\
(3) $\{u^+(\cdot) : u\in  U\}$ and $\{u^-(\cdot) : u\in U\}$ are
equi-left-continuous on $(0, 1]$.
\end{tm}

\begin{re}{\rm Note that $[u]_\al$ is a bounded interval
$[u^-(\al),u^+(\al)]$ for all $u\in E^1$, so condition (3) holds if and only if $\{[u]_\bullet: u\in U\}$
is equi-left-continuous on $(0,1]$. Thus Theorem \ref{fae} is just the case $m=1$ of Theorem \ref{dao}.}
\end{re}

Compare Theorem \ref{cc3} with Theorem \ref{dao}, we know that this is a contradiction because the supremum metric convergence is stronger than the level convergence on $E^m$.

We find that Theorems \ref{dao} and \ref{fae} are incorrect, the following is a counterexample.

\begin{eap}\label{eto}{\rm

Consider a fuzzy number sequence $\{u_n, n=1,2,\ldots\}\subset E^1$ defined by
\begin{gather*}
u_n(\tau)=\left\{\begin{array}{ll}
1,     &  \  \tau=0, \\
\frac{1}{3}+ \frac{2}{3}(1-\tau)^{n},  &  \   0<\tau\leq 1,\\
0, & \ \mbox{otherwise},
\end{array}
\right. \ \  n=1,2,\ldots,
\end{gather*}
then
\begin{gather*}
 [u_n]_\al=\left\{\begin{array}{ll}
[0,1-(\frac{3}{2}\al-\frac{1}{2})^{\frac{1}{n}}],     &  \  \frac{1}{3} <\al \leq 1, \\
\mbox{} [0,1],  &  \ 0\leq \al\leq \frac{1}{3} ,
\end{array}
\right. \ n=1,2,\ldots,
\end{gather*}
for all $\al\in [0,1]$.

We can deduce that
$\{[u_n]_\bullet, n=1,2,\ldots\}$ are equi-left-continuous on $(0,1]$.
In fact, given $\al\in(0,1]$, if $\al\in (\frac{1}{3},1]$, choose a $\delta>0$ such that $\al-\delta>\frac{1}{3}$, then  for all $\beta\in [\al-\delta,\al]$,
\begin{eqnarray}
\lefteqn{H([u_n]_\al, [u_n]_\beta)} \nonumber\\
&=&(\frac{3}{2}\al-\frac{1}{2})^{\frac{1}{n}}-(\frac{3}{2}\beta-\frac{1}{2})^{\frac{1}{n}}\nonumber\\
&\leq& \frac{1}{n} (\frac{3}{2}(\alpha-\delta)-\frac{1}{2})^{\frac{1}{n}-1}(\al-\beta)\nonumber\\
&\leq&  (\frac{3}{2}(\al-\delta)-\frac{1}{2})^{-1}(\al-\beta). \label{dgn}
\end{eqnarray}
If $\al\in [0,\frac{1}{3}]$, then for all $\beta\in [0,\al]$,
\begin{equation}H([u_n]_\al,[u_n]_\beta)=0.\label{dgm}\end{equation}
Combined with \eqref{dgn} and \eqref{dgm}, we know that $\{[u_n]_\bullet, n=1,2,\ldots\}$ are equi-left-continuous on $(0,1]$.

Consider a fuzzy number $u\in E^1$ defined by
\begin{equation*}
u(\tau)=\left\{\begin{array}{ll}
1,  &  \  \tau= 0,\\
\frac{1}{3},     &  \  0 < \tau \leq 1, \\
0, & \ \mbox{otherwise},
\end{array}
\right.\end{equation*}
then \begin{equation*}
u_\al=\left\{\begin{array}{ll}
\{0\},     &  \  \frac{1}{3} <\al \leq 1, \\
\mbox{} [0,1],  &  \ 0\leq \al\leq \frac{1}{3} ,
\end{array}
\right.
\end{equation*}
for all $\al\in [0,1]$. So
\begin{gather}\label{esn}
H(u_\al,[u_n]_\al)=
\left\{\begin{array}{ll}
1-(\frac{3}{2}\al-\frac{1}{2})^{\frac{1}{n}},     &  \  \frac{1}{3} <\al \leq 1, \\
0,  &  \ 0\leq \al\leq \frac{1}{3},
\end{array}
\right.
\end{gather}
and therefore $H(u_\al,[u_n]_\al)\to 0 \; (n\to \infty)$ for all $\al\in [0,1]$ and $d_\infty(u_n,u)=1$ for all $n=1,2,\ldots$, thus we know that $u_n\st{l}{\rightarrow} u$ and $d_\infty(u_n,u)\not\rightarrow 0$. This means that $\{u_n, n=1,2,\ldots \}$ has no limit point in $(E^1, d_\infty)$. So it is a closed set and is not a compact set in $(E^1, d_\infty)$.

Note that $[u_n]_0\subseteq [0,1]$, $n=1,2,\ldots$, i.e. $\{u_n, n=1,2,\ldots\}$ is uniformly support-bounded.
So $\{u_n, n=1,2,\ldots\}$ is a set satisfies conditions (1)-(3) of Theorem \ref{fae}, and it is not a compact set in $(E^1, d_\infty)$. This shows that Theorem \ref{fae} is incorrect.
}
\end{eap}

Theorem 4.1 of \cite{fe} gave a characterization of
compact subsets of all continuous functions from a compact subset
$K$ of a metric space $X$ to $(E^1,d_\infty)$.
However, since it is based on the above theorem, it is wrong too.

\end{document}